\input amstex\documentstyle{amsppt}  
\pagewidth{12.5cm}\pageheight{19cm}\magnification\magstep1
\topmatter
\title A generalization of Steinberg's cross-section \endtitle
\author Xuhua He and George Lusztig\endauthor
\address{Department of Mathematics, Hong Kong University of Science and Technology, Hong Kong and
Department of Mathematics, M.I.T., Cambridge, MA}\endaddress
\thanks{X.H. supported in part by HKRGC grant 601409; G.L. supported in part by the National Science Foundation}
\endthanks
\endtopmatter   
\document

\define\sneq{\subsetneqq}

\define\dz{\dot z}
\define\dw{\dot w}

\define\ds{\dot s}
\define\dy{\dot y}

\define\hW{\hat W}
\define\hx{\hat x}
\define\hy{\hat y}

\define\hw{\hat w}

\define\si{\sim}

\define\sqc{\sqcup}

\define\qua{\quad}

\define\hG{\hat G}

\define\dx{\dot x}

\define\lb{\linebreak}

\define\op{\oplus}

\define\part{\partial}
\define\em{\emptyset}

\define\ra{\rangle}

\define\m{\mapsto}
\define\do{\dots}
\define\la{\langle}
\define\bsl{\backslash}

\define\sub{\subset}    

\define\T{\times}
\define\ti{\tilde}
\define\nl{\newline}
\redefine\i{^{-1}}
\define\fra{\frac}

\define\ot{\otimes}

\define\Aut{\text{\rm Aut}}

\define\Gal{\text{\rm Gal}}

\define\a{\alpha}

\redefine\c{\chi}
\define\g{\gamma}
\redefine\d{\delta}
\define\e{\epsilon}

\define\io{\iota}

\define\p{\pi}
\define\ph{\phi}
\define\ps{\psi}

\define\th{\theta}
\define\k{\kappa}
\redefine\l{\lambda}
\define\z{\zeta}
\define\x{\xi}

\define\Si{\Sigma}
\define\Th{\Theta}

\define\X{\Xi}

\define\bof{\bold f}

\redefine\ss{\bold s}

\define\BB{\bold B}
\define\CC{\bold C}

\define\NN{\bold N}

\define\QQ{\bold Q}

\define\UU{\bold U}

\define\ZZ{\bold Z}

\define\cb{\Cal B}
\define\cc{\Cal C}

\define\co{\Cal O}

\define\car{\Cal R}

\define\ct{\Cal T}

\define\cx{\Cal X}

\define\fB{\frak B}

\define\tn{\ti n}

\define\tH{\ti H}

\define\tU{\ti U}

\define\tX{\ti X}

\define\dUU{\dot{\UU}}
\define\hUU{\hat{\UU}}
\define\dBB{\dot{\BB}}
\define\bco{\bar{\co}}
\define\tfB{\ti{\fB}}
\define\AX{Ax}
\define\BR{BR}
\define\DL{DL}
\define\GP{GP}
\define\GKP{GKP}
\define\GM{GM}
\define\GRR{G1}
\define\GR{G2}
\define\HE{He}
\define\QG{L1}
\define\LU{L2}
\define\WEU{L3}
\define\HO{L4}
\define\XW{L5}
\define\SE{Se}
\define\ST{St}

\head Introduction \endhead
\subhead 0.1\endsubhead
Let $G$ be a connected semisimple algebraic group over an algebraically closed field. Let $B,B^-$ be two opposed
Borel subgroups of $G$ with unipotent radicals $U,U^-$ and let $T=B\cap B^-$, a maximal torus of $G$. Let $NT$ be
the normalizer of $T$ in $G$ and let $W=NT/T$ be the Weyl group of $T$, a finite Coxeter group with length 
function $l$. For $w\in W$ let $\dw$ be a representative of $w$ in $NT$. The following result is due to Steinberg
\cite{\ST, 8.9} (but the proof in {\it loc.cit.} is omitted): if $w$ is a Coxeter element of minimal length in 
$W$ then (i) the conjugation action of $U$ on $U\dw U$ has trivial isotropy groups and (ii) the subset 
$(U\cap\dw U^-\dw\i)\dw$ meets any $U$-orbit on $U\dw U$ in exactly one point; in particular, (iii) the set of 
$U$-orbits on $U\dw U$ is naturally an affine space of dimension $l(w)$.

More generally, assuming that $w$ is any elliptic element of $W$ of minimal length in its conjugacy class, it is 
shown in \cite{\WEU} that (i) holds and, assuming in addition that $G$ is of classical type, it is shown in
\cite{\XW} that (iii) holds. In this paper we show for any $w$ as above and any $G$ that (ii) (and hence (iii)) 
hold, see Theorem 3.6(ii) (actually we take $\dw$ of a special form but then the result holds in general since 
any representative of $w$ in $NT$ is of the form $t\dw t\i$ for some $t\in T$). We also prove analogous 
statements in some twisted cases, involving an automorphism of the root system or a Frobenius map (see 3.6) 
and a version over $\ZZ$ of these statements using the results in \cite{\LU} on groups over $\ZZ$. Note that the 
proof of (ii) given in this paper uses (as does the proof of (i) in \cite{\WEU}) a result in \cite{\GP, 3.2.7} 
and a weak form of the existence of "good elements" \cite{\GM} in an elliptic conjugacy class in $W$ which (for 
exceptional types) rely on a computer.

\subhead 0.2\endsubhead
Let $w$ be an elliptic element of $W$ which has minimal length in its conjugacy class $C$ and let $\g$ be the
unipotent class of $G$ attached to $C$ in \cite{\WEU}. Recall that $\g$ has codimension $l(w)$ in $G$. As an 
application of our results we construct (see 4.2(a)) a closed subvariety $\Si$ of $G$ isomorphic to an affine 
space of dimension $l(w)$ such that $\Si\cap\g$ is a finite set with a transitive action of a certain finite group
whose order is divisible only by the bad primes of $G$. In the case where $C$ is the Coxeter class, $\Si$ reduces 
to Steinberg's cross section \cite{\ST} which intersects the regular unipotent class in $G$ in exactly one 
element. 

\subhead 0.3\endsubhead
Recently, A. Sevostyanov \cite{\SE} proved statements similar to (i),(ii),(iii) in 0.1 for a certain type of Weyl 
group elements assuming that the ground field is $\CC$. It is not clear to us what is the relation of the Weyl 
group elements considered in \cite{\SE} with those considered in this paper. 

\subhead 0.4\endsubhead
The following (unpublished) example of N. Spaltenstein, dating from the late 1970's, shows that the statement (i)
(for Coxeter elements) in 0.1 can be false if the assumption of minimal length is dropped: the elements 
$$\align&\dw=\left(\matrix 
0&0&1&0&0&0\\
0&0&0&0&1&0\\
0&1&0&0&0&0\\
0&0&0&0&0&1\\
0&0&0&1&0&0\\
1&0&0&0&0&0\endmatrix\right),\qua u_x=\left(\matrix 
1&0&0&0&0&0\\
0&1&0&x&0&x\\
0&0&1&x&0&x\\
0&0&0&1&0&0\\
0&0&0&0&1&0\\
0&0&0&0&0&1\endmatrix\right),\\&
y=\left(\matrix 
1&0&-1&0&0&0\\
0&1&1&0&0&0\\
0&0&1&0&0&0\\
0&0&0&1&1&-1\\
0&0&0&0&1&0\\
0&0&0&0&0&1\endmatrix\right)\endalign$$
of $GL_6(\CC)$ (with $x\in\CC$) satisfy $u_xy\dw u_x\i=y\dw$; hence if $U$ is the group of upper triangular 
matrices in $GL_6(\CC)$ then in the conjugation action of $U$ on $U\dw U$, the isotropy group of $y\dw$ contains 
the one parameter group $\{u_x;x\in\CC\}$.

\head 1. Polynomial maps of an affine space to itself\endhead
\subhead 1.1\endsubhead
Let $\cc$ be the class of commutative rings with $1$. Let $N$ be an integer $\ge1$.
A family $(f_A)_{A\in\cc}$ of maps $f_A:A^N@>>>A^N$ is said to be {\it polynomial} if there exist (necessarily 
unique) polynomials with integer coefficients $$f_1(X_1,\do,X_N),\do,f_N(X_1,\do,X_N)$$ in the indeterminates 
$X_1,\do,X_N$ such that for any $A\in\cc$, $f_A$ is the map 
$$(a_1,\do,a_N)\m(f_1(a_1,\do,a_N),\do,f_N(a_1,\do,a_N)).$$

For such a family we define for any $A\in\cc$ an $A$-algebra homomorphism 
$f^*_A:A[X_1,\do,X_N]@>>>A[X_1,\do,X_N]$ by $f^*_A(X_i)=f_i(X_1,\do,X_N)$ for all $i\in[1,N]$ (here we view 
$f_i$ as an element of $A[X_1,\do,X_N]$ using the obvious ring homomorphism $\ZZ@>>>A$). We have the following 
result.

\proclaim{Proposition 1.2} Assume that $f_A:A^N@>>>A^N$ ($A\in\cc$) is a polynomial family such that $f_A$ is
injective for any $A\in\cc$. Then:

(i) $f_A$ is bijective in the following cases: (a) $A$ is finite; (b) $A=\ZZ_p$, the ring of $p$-adic integers 
($p$ is a prime number); (c) $A$ is a perfect field; (d) $A$ is the ring of rational numbers which have no $p$ in
the denominator ($p$ is a prime number); (e) $A=\ZZ$.

(ii) If $A$ is an algebraically closed field then $f_A$ is an isomorphism of algebraic varieties.
\endproclaim
We prove (i). In case (a), $A^N$ is a finite set and the result follows.

Assume that $A$ is as in (b). For any $s\ge1$ let $A_s=\ZZ/p^s\ZZ$ (a finite ring) and let $l_s:A@>>>A_s$ be the 
obvious homomorphism. Let $\x\in A^N$. Let $\x_s=l_s^N(\x)\in A_s^N$. Using (a) for $A_s$ we set
$\x'_s=f_{A_s}\i(\x_s)\in A_s^N$. Let $\x'\in A^N$ be the unique element such that for any $s\ge1$, the image of 
$\x'$ under $l_s^N:A^N@>>>A^N_s$ is equal to $\x'_s$. Let $\ti\x=f_A(\x')$. Then for any $s\ge1$, $\x,\ti\x$ have
the same image under $l_s^N:A^N@>>>A_s^N$. Hence $\ti\x=\x$. Thus $f_A$ is surjective, as desired.

Assume that $A$ is as in (c). Let $A'$ be an algebraic closure of $A$. By \cite{\BR} (see also \cite{\AX}, 
\cite{\GRR, 10.4.11}), $f_{A'}$ is bijective. Let $\x\in A^N$. Since $A^N\sub A'{}^N$, we can view $\x$ as an 
element of $A'{}^N$ so that $\x'=f_{A'}\i(\x)\in A'{}^N$ is defined. For any $\g\in\Gal(A'/A)$ we have 

$f_{A'}(\g(\x'))=\g(f_{A'}(\x'))=\g(\x)=\x=f_{A'}(\x')$.
\nl
(The obvious action of $\g$ on $A'{}^N$ is denoted again by $\g$.) Using the injectivity of $f_{A'}$ we deduce 
that $\g(\x')=\x'$. Since this holds for any $\g$ and $A$ is perfect, it follows that $\x'\in A^N$. We have 
$f_A(\x')=\x$. Thus $f_A$ is surjective, as desired.

Assume that $A$ is as in (d). Let $A_0=\QQ_p$, the field of $p$-adic numbers. We can view $A$ as the intersection
of two subrings of $A_0$, namely $A_1=\QQ$ and $A_2=\ZZ_p$. Now $f_{A_i}$ is bijective for $i=0,1,2$ by (b),(c). 
Let $\x\in A^N$. For $i=0,1,2$ we set $\x'_i=f_{A_i}\i(\x)\in A_i^N$. Clearly, $\x'_1=\x'_0=\x'_2$ hence
$\x'_0\in A^N$ and $f_A(\x'_0)=\x$. Thus $f_A$ is surjective, as desired.

Assume that $A=\ZZ$. We can view $A$ as $\cap_pA_p$ where $p$ runs over the set of prime numbers and $A_p$ is the
ring denoted by $A$ in (d); the intersection is taken in the field $\QQ$. Now $f_{A_p}$ is bijective for any $p$ 
(see (d)) and $f_\QQ$ is bijective by (c). Let $\x\in A^N$. We set $\x'_p=f_{A_p}\i(\x)\in A_p^N$ and 
$\x'=f_\QQ\i(\x)\in\QQ^N$. Clearly, $\x'_p=\x'$ for all $p$. Hence $\x'\in A^N$ and $f_A(\x')=\x$. Thus $f_A$ is 
surjective, as desired. This proves (i).

Now assume that $A$ is as in (ii). Since $f_{A_1}:A_1^N@>>>A_1^N$ is bijective for any algebraically closed field
$A_1$ (see (i)), it is enough to show that the morphism $f_A$ is \'etale (see \cite{\GR, 17.9.1}). Let 
$A'=A\op A$ regarded as an $A$-algebra with multiplication $(a,b)(a',b)=(ab,ab'+a'b)$. The unit element is 
$1=(1,0)$. We set $T=(0,1)$. Then $(a,b)=a+bT$ and $T^2=0$. Then $f_{A'}$ is defined.
There exist polynomials $f_k(X_1,\do,X_N)$ ($k\in[1,N]$) with coefficients in $\ZZ$ such that
$$\align&f_{A'}(a_1+b_1T,\do,a_N+b_NT)\\&=(f_1(a_1+b_1T,\do,a_N+b_NT),\do,f_N(a_1+T,\do,a_N+Tb_N))\\&
=(f_1(a_*)+\sum_{k=1}^N\fra{\partial f_1}{\partial X_k}(a_*)b_kT,\do,f_N(a_*)+\sum_{k=1}^N
\fra{\partial f_N}{\partial X_k}(a_*)b_kT)\endalign$$
for any $a_*=(a_1,\do,a_N)\in A^N$, $b_*=(b_1,\do,b_N)\in A^N$. Since $f_{A'}$ is injective,
we see that for any $a_*\in A^N$, the (linear) map $A^N@>>>A^N$ given by
$$b_*\m(\sum_{k=1}^N\fra{\partial f_1}{\partial X_k}(a_*)b_k,\do,
\sum_{k=1}^N\fra{\partial f_N}{\partial X_k}(a_*)b_k)$$
is injective. It follows that for any $a_*\in A^N$, the $N\T N$-matrix with $(j,k)$-entries 
$\fra{\partial f_j}{\partial X_k}(a_*)$ is nonsingular. This shows that $f_A$ is \'etale. This proves (ii).
The proposition is proved.

\proclaim{Proposition 1.3} Assume that $f_A:A^N@>>>A^N$ ($A\in\cc$) and $f'_A:A^N@>>>A^N$ ($A\in\cc$) are two
polynomial families such that $f'_Af_A=1$ for all $A\in\cc$. Then for any $A\in\cc$, $f^*_A$ is an $A$-algebra 
isomorphism and $f_A:A^N@>>>A^N$ is bijective.
\endproclaim
Since $f'_Af_A=1$ we have $f_A^*f'_A{}^*=1$ and $f_A$ is injective for any $A$.
Using 1.2 we see that $f_\ZZ$ is bijective hence $f_\ZZ f'_\ZZ=1$. Let $\x_A=f_Af'_A$. Then $(\x_A)_{A\in\cc}$ 
is a polynomial family and $\x_\ZZ=1$. Thus $\x^*_\ZZ(X_i)=X_i$ for any $i$ (there is at most one element of 
$\ZZ[X_1,\do,X_N]$ with prescribed values at any $(x_1,\do,x_N)\in\ZZ^N$). We see that the polynomials with 
integer coefficients which define $\x$ are $X_1,\do,X_N$. It follows that $\xi_A=1$ for any $A\in\cc$; hence 
$f_Af'_A=1$ and $f_A$ is a bijection. Also, since $\x^*_A=1$ we have $f'_A{}^*f_A^*=1$ for any $A$ hence $f^*_A$ 
is an isomorphism. The proposition is proved.

\head 2. Reductive groups over a ring\endhead
\subhead 2.1\endsubhead
We fix a root datum $\car$ as in \cite{\QG, 2.2}. This consists of two free abelian groups of finite type $Y,X$ 
with a given perfect pairing $\la,\ra:Y\T X\to\ZZ$ and a finite set $I$ with given imbeddings $I@>>>Y\qua(i\m i)$
and $I@>>>X\qua(i\m i')$ such that $\la i,i'\ra=2$ for all $i\in I$ and $\la i,j'\ra\in-\NN$ for all $i\ne j$ in 
$I$; in addition, we are given a symmetric bilinear form $\ZZ[I]\T\ZZ[I]\to\ZZ$, $\nu,\nu'\m\nu\cdot\nu'$ such 
that $i\cdot i\in2\ZZ_{>0}$ for all $i\in I$ and $\la i,j'\ra=2i\cdot j/i\cdot i$ for all $i\ne j$ in $I$. 
We assume that the matrix $M=(i\cdot j)_{i,j\in I}$ is positive definite.

Let $W$ be the (finite) subgroup of $\Aut(X)$ generated by the involutions $s_i:x-\la i,x\ra i'\qua (i\in I)$. 
For $i\ne j$ in $I$ let $n_{i,j}=n_{j,i}$ be the order of $s_is_j$ in $W$. Note that $W$ is a (finite) Coxeter 
group with generators $S:=\{s_i;i\in I\}$; let $l:W@>>>\NN$ be the standard length function. Let $w_I$ be the 
unique element of maximal length of $W$. For $J\sub I$ let $W_J$ be the subgroup of $W$ generated by
$\{s_i;i\in J\}$. Let $\cx$ be the set of all sequences $i_1,i_2,\do,i_n$ in $I$ such that 
$s_{i_1}s_{i_2}\do s_{i_n}=w_I$ and $l(w_I)=n$.

\subhead 2.2\endsubhead
Now (until the end of 2.10) we fix $A\in\cc$. Let $\dUU_A$ be the $A$-algebra attached to $\car$ and to $A$ (with 
$v=1$) in \cite{\QG, 31.1.1}, where it is denoted by ${}_A\dUU$. As in {\it loc.cit.} we denote the canonical 
basis of the $A$-module $\dUU_A$ by $\dBB$. For $a,b,c\in\dBB$ we define $m_{a,b}^c\in A$, $\hat m^{a,b}_c\in A$ 
as in \cite{\LU, 1.5}; in particular we have $ab=\sum_{c\in\dBB}m_{a,b}^cc$. Let $\hUU_A$ be the $A$-module 
consisting of all formal linear combinations $\sum_{a\in\dBB}n_aa$ with $n_a\in A$. There is a well defined 
$A$-algebra structure on $\hUU_A$ such that 
$$(\sum_{a\in\dBB}n_aa)(\sum_{b\in\dBB}\tn_bb)=\sum_{c\in\dBB}r_cc$$
where $r_c=\sum_{(a,b)\in\dBB\T\dBB}m_{a,b}^cn_a\tn_b$. (See \cite{\LU, 1.11}.) It has unit element
$1=\sum_{\z\in X}1_\z$ where for any $\z\in X$ the element $1_\z\in\dBB$ is defined as in \cite{\QG, 31.1.1}. Let
$\e:\hUU_A@>>>A$ be the algebra homomorphism given by $\sum_{a\in\dBB}n_aa\m n_{1_0}$. We identify $\dUU_A$ with 
the subalgebra of $\hUU_A$ consisting of finite $A$-linear combinations of elements in $\dBB$.

Let $\bof_A$ be the $A$-algebra with $1$ associated to $M$ and $A$ (with $v=1$) in \cite{\QG, 31.1.1}, where it 
is denoted by ${}_A\bof$. Let $\BB$ be the canonical basis of the $A$-module $\bof_A$ (see \cite{\QG, 31.1.1}). 
For $i\in I,c\in\NN$, the element $\th_i^{(c)}$ (see \cite{\QG, 1.4.1, 31.1.1}) of $\bof_A$ is contained in 
$\BB$. Let $(x\ot x'):u\m x^-ux'{}^+$ be the $\bof_A\ot_A\bof_A^{opp}$-module structure on $\dUU_A$ considered 
in \cite{\QG, 31.1.2}; we write $ux'{}^+$ instead of $1^-ux'{}^+$ and $x^-u$ instead of $x^-u1^+$.

The elements $\{1_\z b^+;b\in\BB,\z\in X\}$ (resp. $\{b^-1_\z;b\in\BB,\z\in X\}$) of $\dUU_A$ are distinct and 
form a subset of $\dBB$. Let $\hUU^+_A$ (resp. $\hUU^-_A$) be the $A$-submodule of $\hUU_A$ consisting of elements
$$\sum_{b\in\BB,\z\in X}n_b(1_\z b^+)\text{(resp. }\sum_{b\in\BB,\z\in X}n_b(b^-1_\z))$$ 
with $n_b\in A$. 

\subhead 2.3\endsubhead
As in \cite{\LU, 4.1}, let $G_A$ be the set of all $\sum_{a\in\dBB}n_aa\in\hUU_A$ such that $n_{1_0}=1$ and
$\sum_{c\in\dBB}\hat m^{a,b}_cn_c=n_an_b$ for all $a,b\in\dBB$ (the last sum is finite by \cite{\LU, 1.16}). Note 
that $G_A$ is a subgroup of the group of invertible elements of the algebra $\hUU_A$. As in \cite{\LU, 4.2} we 
set $U_A=G_A\cap\hUU_A^+$, $U_A^-=G_A\cap\hUU_A^-$ (in {\it loc.cit.} $U_A,U_A^-$ are denoted by 
$G_A^{>0},G_A^{<0}$). Let $T_A$ the set of elements of $G_A$ of the form $\sum_{\l\in X}n_\l1_\l$ with 
$n_\l\in A$. As shown in {\it loc.cit.}, $U_A,U_A^-,T_A$ are subgroups of $G_A$. We note that

(a) {\it multiplication in $G_A$ defines an injective map $U_A^-\T T_A\T U_A@>>>G_A$.}
\nl
This statement appears in \cite{\LU, 4.2(a)} but the line

"$\hUU^-_A\cap\hUU_A^{>0}=\{1\}$. Thus, $\x'_3\x_3\i=1$ so that $\x'_3=\x_3$." 
\nl
in the proof in {\it loc. cit.} should be replaced by:

"$\hUU^-_A\cap\hUU_A^{>0}=\{a1;a\in A\}$. Thus, $\x'_3\x_3\i=a1$ for some $a\in A$. Since 
$\e(\x'_3\x_3\i)=1$, we have $a=1$ so that $\x'_3=\x_3$."
\nl
From (a) we deduce that

(b) {\it $U_A\cap U_A^-=\{1\}$.}

\subhead 2.4\endsubhead
For any $i\in I,h\in A$ we set
$$x_i(h)=\sum_{c\in\NN,\l\in X}h^c1_\l\th_i^{(c)+}\in\hUU_A,$$
$$y_i(h)=\sum_{c\in\NN,\l\in X}h^c\th_i^{(c)-}1_\l\in\hUU_A.$$
By \cite{\LU, 1.18(a)} we have $x_i(h)\in U_A,y_i(h)\in U_A^-$. For any $i\in I$ we set
$$\ds_i=\sum_{a\in\NN,b\in\NN,\l\in X;\la i,\l\ra=a+b}(-1)^a\th_i^{(a)-}1_\l\th_i^{(b)+}\in\hUU_A.$$
(Note that $\th_i^{(a)-}1_\l\th_i^{(b)+}\in\dBB$ whenever $a+b=\la i,\l\ra$.) By \cite{\LU, 2.2(e)} we have 
$\ds_i\in G_A$. (In {\it loc.cit.}, $\ds_i$ is denoted by $s''_{i,1}$.) From \cite{\LU, 2.4} we see that if 
$i,j\in I$, $i\ne j$, then

(a) $\ds_i\ds_j\ds_i\do=\ds_j\ds_i\ds_j\do$  in $G_A$ (both products have $n_{i,j}$ factors).
\nl
For $i\in I$ we have

(b) $\ds_i^2=t_i(-1)$
\nl
where $t_i(-1)=\sum_{\l\in X}(-1)^{\la i,\l\ra}1_\l\in T_A$ normalizes $U_A,U_A^-$. (See 
\cite{\LU, 2.3(b), 4.4(a), 4.3(a)}.)

For $i\in I,h\in A$ we have 

(c) $\ds_i\i x_i(h)\ds_i=y_i(-h)$.
\nl
This follows from the definitions using \cite{\LU, 2.3(c)}.

For any $w\in W$ we set $\dw=\ds_{i_1}\ds_{i_2}\do\ds_{i_k}\in G_A$ where $i_1,\do,i_k$ in $I$ are such that 
$s_{i_1}s_{i_2}\do s_{i_k}=w$, $k=l(w)$. This is well defined, by (a).

\subhead 2.5\endsubhead
We have the following result (see \cite{\LU, 4.7(a)}):

(a) {\it Let $w\in W$ and let $i\in I$ be such that $l(ws_i)=l(w)+1$. Let $h\in A$. We have 
$\dw x_i(h)\dw\i\in U_A$.}
\nl
The proof of the following result is similar to that of (a):

(b) {\it Let $z\in W$ and let $i\in I$ be such that $l(s_iz)=l(z)+1$. Let $h\in A$. Then we have 
$\dz\i y_i(h)\dz\in U_A^-$.}

\subhead 2.6\endsubhead
Let $(i_1,\do,i_n)\in\cx$. For any $(h_1,h_2,\do,h_n)\in A^n$ we have
$$\align&x_{i_1}(h_1)\ds_{i_1}x_{i_2}(h_2)\ds_{i_2}\do x_{i_n}(h_n)\ds_{i_n}\dw_I\i\\&=
x_{i_1}(h_1)(\ds_{i_1}x_{i_2}(h_2)\ds_{i_1}\i)\do
(\ds_{i_1}\ds_{i_2}\do\ds_{i_{n-1}}x_{i_n}(h_n)\ds_{i_{n-1}}\i\do\ds_{i_2}\i\ds_{i_1}\i)\in U_A.\endalign$$
(We use 2.5(a).) From \cite{\LU, 4.8(a)} we see that 

(a) {\it the map $A^n@>>>U_A$, $(h_1,h_2,\do,h_n)\m 
x_{i_1}(h_1)\ds_{i_1}x_{i_2}(h_2)\ds_{i_2}\do x_{i_n}(h_n)\ds_{i_n}\dw_I\i$ is a bijection.}

\subhead 2.7\endsubhead
Let $w\in W$. Let 
$$U^w_A=U_A\cap\dw U_A^-\dw\i,\qua {}^wU_A=U_A\cap\dw U_A\dw\i;$$
these are subgroups of $U_A$. We can find $(i_1,\do,i_n)\in\cx$ such that $s_{i_1}s_{i_2}\do s_{i_k}=w$ where 
$k=l(w)$. Let $(h_1,h_2,\do,h_n)\in A^n$ and let $u\in U_A$ be its image under the map 2.6(a). We have $u=u'u''$ 
where $u'=r_1r_2\do r_k$, $u''=r_{k+1}r_{k+2}\do r_n$ and
$$r_m=\ds_{i_1}\ds_{i_2}\do\ds_{i_{m-1}}x_{i_m}(h_m)\ds_{i_{m-1}}\i\do\ds_{i_2}\i\ds_{i_1}\i\in U_A$$
for $m\in[1,n]$. If $m\in[1,k]$ we have 
$$\align&\dw\i r_m\dw=\ds_{i_k}\i\ds_{i_{k-1}}\i\do\ds_{i_m}\i x_{i_m}(h)\ds_{i_m}\do\ds_{i_{k-1}}\ds_{i_k}\\&=
\ds_{i_k}\i\ds_{i_{k-1}}\i\do\ds_{i_{m+1}}\i y_{i_m}(-h)\ds_{i_{m+1}}\do\ds_{i_{k-1}}\ds_{i_k}\in U_A^-\endalign$$
(we have used 2.4(c), 2.5(b)). Hence $\dw\i u'\dw\in U_A^-$ and $u'\in U^w_A$. If $m\in[k+1,n]$ we have 
$$\dw\i r_m\dw=\ds_{i_{k+1}}\ds_{i_{k+2}}\do\ds_{i_{m-1}}x_{i_m}(h_m)\ds_{i_{m-1}}\i\do\ds_{i_{k+2}}\i
\ds_{i_{k+1}}\i\in U_A$$
(we have used 2.5(a)). Hence $\dw\i u''\dw\in U_A$ and $u''\in {}^wU_A$.
If we assume that $u\in U^w_A$ then $u''=u'{}\i u\in U^w_A$ hence $\dw\i u''\dw\in U_A^-$. We have also
$\dw\i u''\dw\in U_A$. Since $U_A\cap U_A^-=\{1\}$ (by 2.3(b))) we have $\dw\i u''\dw=1$ and $u''=1$, $u=u'$.
Conversely, if $u=u'$ then, as we have seen, we have $u\in U^w_A$. 

If we assume that $u\in{}^wU_A$ then $u'=uu''{}\i\in{}^wU_A$ hence $\dw\i u'\dw\in U_A$. We have also
$\dw\i u'\dw\in U_A^-$. Since $U_A\cap U_A^-=\{1\}$ (by 2.3(b)) we have $\dw\i u'\dw=1$ and $u'=1$, $u=u''$.
Conversely, if $u=u''$ then, as we have seen, we have $u\in{}^wU_A$. Thus we have the following results:

(a) {\it the restriction of the map 2.6(a) to $A^k$ (identified with 
$$\{(h_1,h_2,\do,h_n)\in A^n;h_{k+1}=h_{k+2}=\do=h_n=0\}$$
 in the obvious way) defines a bijection $A^k@>\si>>U_A^w$;}

(b) {\it the restriction of the map 2.6(a) to $A^{n-k}$ (identified with 
$$\{(h_1,h_2,\do,h_n)\in A^n;h_1=h_2=\do=h_k=0\}$$
 in the obvious way) defines a bijection $A^{n-k}@>\si>>{}^wU_A$.}
\nl
Using (a),(b) and 2.6 we see also that

(c) {\it multiplication defines a bijection $U^w_A\T{}^wU_A@>\si>>U_A$.} 
\nl
We show:

(d) {\it multiplication in $G_A$ defines a bijection $(U^w_A\dw)\T U_A@>\si>>U_A\dw U_A$.} 
\nl
Using (c) we see that $U_A\dw U_A=U^w_A({}^wU_A)\dw U_A=U^w_A\dw(\dw\i{}^wU_A\dw) U_A\sub U^w_A\dw U_A$. Thus the
map in (d) is surjective. Assume now that $u_1,u_2\in U^w_A$ and $u'_1,u'_2\in U_A$ satisfy 
$u_1\dw u'_1=u_2\dw u'_2$. Then $\dw\i u_2\i u_1\dw=u'_2u'_1{}\i$ is both in $U_A^-$ and in $U_A$ hence by 2.3(b)
it is $1$. Thus $u_1=u_2$ and $u'_1=u'_2$. Thus the map in (d) is injective. This proves (d).

Combining (d) with (a) and 2.6(a), we obtain a bijection

(e) $A^k\T A^n@>\si>>U_A\dw U_A$, 
$$\align&((h_1,\do,h_k),(h'_1,\do,h'_n))\m(x_{i_1}(h_1)\ds_{i_1}x_{i_2}(h_2)\ds_{i_2}\do x_{i_k}(h_k)
\ds_{i_k}\dw\i)\dw\\&(x_{i_1}(h'_1)\ds_{i_1}x_{i_2}(h'_2)\ds_{i_2}\do x_{i_n}(h'_n)\ds_{i_n}\dw_I\i).\endalign$$
We can reformulate (a) as follows.

(f) {\it If $j_1,j_2,\do,j_k$ is a sequence in $I$ such that $w=s_{j_1}\do s_{j_k}$, $l(w)=k$, then the map 
$A^k@>>>U^w_A\dw$, 
$$(h_1,h_2,\do,h_k)\m x_{i_1}(h_1)\ds_{i_1}x_{i_2}(h_2)\ds_{i_2}\do x_{i_k}(h_k)\ds_{i_k}$$
 is a bijection.}
\nl
For any sequence $w_*=(w_1,w_2,\do,w_r)$ in $W$ we set
$$U(w_*)=(U^{w_1}_A\dw_1)\T(U^{w_2}_A\dw_2)\T\do\T(U^{w_r}_A\dw_r),$$
$$\dot U(w_*)=(U_A\dw_1U_A)\T(U_A\dw_2U_A)\T\do\T(U_A\dw_rU_A).$$
We have an obvious inclusion $U(w_*)\sub\dot U(w_*)$. Let $\tU(w_*)$ be the set of orbits of the 
$U_A^{r-1}$-action 
$$(u_1,u_2,\do,u_{r-1}):(g_1,g_2,\do,g_r)\m(g_1u_1\i,u_1g_2u_2\i,\do,u_{r-1}g_r)$$
on $\dot U(w_*)$. Let $\k_{w_*}:\dot U(w_*)@>>>\tU(w_*)$ be the obvious surjective map; for \lb
$(g_1,\do,g_r)\in\dot U(w_*)$ we set $[g_1,\do,g_r]=\k_{w_*}(g_1,\do,g_r)$.

The following result is an immediate consequence of (f).

(g) {\it Let $w\in W$ and let $w_*=(w_1,w_2,\do,w_r)$ be a sequence in $W$ such that $w=w_1w_2\do w_r$, 
$l(w)=l(w_1)+l(w_2)+\do+l(w_r)$. Then multiplication in $G_A$ defines a bijection
$\ph_{w_*}:U(w_*)@>\si>>U^w_A\dw$.}
\nl
We show:

(h) {\it Let $x,y\in W$ be such that $l(xy)=l(x)+l(y)$. Let $x_*=(x,y)$. Then multiplication in $G_A$ defines
a bijection $\tU(x_*)@>\si>>U_A\dx\dy U_A$.}
\nl
This follows from (d) and (g) since
$$\tU(x_*)\cong(U^x_A\dx)\T(U_A\dy U_A)=(U^x_A\dx)\T(U^y_A\dy U_A)=U^{xy}_A\dx\dy U_A.$$
Using (h) repeatedly we obtain:

(i) {\it In the setup of (g), multiplication in $G_A$ defines a bijection
$\ps_{w_*}:\tU(w_*)@>\si>>U_A\dw U_A$. This bijection is compatible with the $U_A\T U_A$ actions 

$(u,u'):[g_1,g_2,\do,g_r]\m[ug_1,g_2,\do,g_ru'{}\i]$
\nl
on $\tU(w_*)$ and $(u,u'):g\m ugu'{}\i$ on $U_A\dx U_A$.}

\subhead 2.8\endsubhead
Let $y_*=(y_1,y_2,\do,y_t)$ ($t\ge2$) be a sequence in $W$. Let $\x\in\tU(y_*)$. We show:

(a) {\it for any $s\in[1,t-1]$ there exists a representative of $\x$ in
$$(U_A^{y_1}\dy_1)\T\do\T(U_A^{y_{s-1}}\dy_{s-1})\T(U_A\dy_sU_A)\T\do\T(U_A\dy_tU_A).$$}
\nl
We argue by induction on $s$. Let $(g_1,g_2,\do,g_t)\in\dot U(y_*)$ be a representative of $\x$. We have
$g_1=u\dy_1u'$ where $u\in U_A^{y_1},u'\in U_A$ (see 2.7(d)). Then $(u\dy_1,u'g_2,g_3,\do,g_t)$ 
represents $\x$ and is as in (a) with $s=1$. Now assume that $s\ge2$ and that $(g_1,g_2,\do,g_t)$ is a 
representative of $\x$ as in (a) with $s$ replaced by $s-1$. We have 
$g_{s-1}=u\dy_{s-1}u'$ with 
$u\in U_A^{y_{s-1}},u'\in U_A$ (see 2.7(d)). Then 
$$(g_1,g_2,\do,g_{s-2},u\dy_{s-1},u'g_s,g_{s+1},\do,g_t)$$
is a representative of $\x$ as in (a). This completes the inductive proof of (a).

We show:

(b) {\it there exist unique elements $u_1\in U_A^{y_1},\do,u_t\in U_A^{y_t}$ and $u\in U_A$ such that
$\x=[u_1\dy_1,\do,u_{t-1}\dy_{t-1},u_t\dy_tu]$.}
\nl
The existence of these elements folows from (a) with $s=t-1$ and 2.7(d). We prove uniqueness. Assume that
$u_1,u'_1\in U_A^{y_1},\do,u_t,u'_t\in U_A^{y_t}$ and $u,u'\in U_A$ are such that
$$[u_1\dy_1,\do,u_{t-1}\dy_{t-1},u_t\dy_tu]=[u'_1\dy_1,\do,u'_{t-1}\dy_{t-1},u'_t\dy_tu']=\x.$$
Then there exist $v_1,v_2,\do,v_{t-1}$ in $U_A$ such that
$$\align&u'_1\dy_1=u_1\dy_1v_1, u'_2\dy_2=v_1\i u_2\dy_2v_2, \do,
u'_{t-1}\dy_{t-1}=v_{t-2}\i u_{t-1}\dy_{t-1}v_{t-1},\\& 
u'_t\dy_tu'=v_{t-1}\i u_t\dy_tu.\endalign$$
The first of these equations implies (using 2.7(d)) that $v_1=1$ and $u'_1=u_1$. Then the second equation becomes
$u'_2\dy_2=u_2\dy_2v_2$; using again 2.7(d) we deduce that $v_2=1$ and $u'_2=u_2$. Continuing in this way
we get $v_1=\do=v_{t-1}=1$ and $u_i=u'_i$ for $i\in[1,t-1]$. We then have $u'_t\dy_tu'=u_t\dy_tu$. Using 
2.7(d) we deduce $u'_t=u_t$, $u'=u$. This proves (b).

We now define $\z:\tU(y_*)@>>>U_A$ by $\z(\x)=u$ where $u$ is as in (b).

\subhead 2.9\endsubhead
Let $\hW$ be the braid monoid attached to $W$ and let $w\m\hw$ be the canonical map $W@>>>\hW$. Let 
$x_*=(x_1,\do,x_s)$, $y_*=(y_1,\do,y_t)$ be two sequences in $W$ such that $\hx_1\do\hx_s=\hy_1\do\hy_t$ in 
$\hW$. We show:

(a) {\it there exist bijections $H:U(x_*)@>\si>>U(y_*)$, $\tH:\tU(x_*)@>\si>>\tU(y_*)$ such that
$\k_{y_*}H=\tH\k_{x_*}$ and such that $\tH$ is compatible with the $U_A\T U_A$-actions (as in 2.7(i)).}
\nl
Applying 2.7(g) and 2.7(i) with $w$ replaced by $x_a$ (resp. $y_b$) and $w_1,w_2,\do,w_r$ replaced by a sequence 
of simple reflections whose product is a reduced expression of $x_a$ (resp. $y_b$) we see that the general case 
is reduced to the case where $l(x_1)=\do=l(x_s)=l(y_1)=\do=l(y_t)=1$; in this case we must have $s=t$. Since
$\hx_1\do\hx_s=\hy_1\do\hy_s$ we can find a sequence $\ss^1_*,\ss^2_*,\do,\ss^m_*$ where each $\ss^p_*$ is a 
sequence in $S$, $\ss^1_*=x_*$, $\ss^m_*=y_*$ and for any $r\in[1,m-1]$, $\ss^{r+1}_*$ is related to $\ss^r_*$ as
follows:

(b) for some $i\ne j$ in $I$ and some $e$ such that $[e+1,e+u]\sub[1,s]$ ($n_{ij}=u$) we have 

$\ss^r_c=\ss^{r+1}_c$ if $c\in[1,s]-[e+1,e+u]$,

$(\ss^r_{e+1},\ss^r_{e+2},\do,\ss^r_{e+u})=s_*:=(s_i,s_j,s_i,\do)$,

$(\ss^{r+1}_{e+1},\ss^{r+1}_{e+2},\do,\ss^{r+1}_{e+u})=s'_*:=(s_j,s_i,s_j,\do)$.
\nl
Note that each pair $\ss^r_*,\ss^{r+1}_*$ satisfies the same assumptions as the pair $x_*,y_*$. If (a) can be 
proved for each $\ss^r,\ss^{r+1}$ then (a) would follow also for $x_*,y_*$ (by taking appropriate compositions of
maps like $H$ or maps like $\tH$). It is then enough to prove (a) assuming that $x_*=\ss^r_*$, $y_*=\ss^{r+1}_*$ 
are as in (b). Let $w=s_is_js_i\do=s_js_is_j\do$ ($u$ factors). Let 
$$x'_*=(x_1,\do,x_e)=(y_1,\do,y_e),\qua x''_*=(x_{e+u+1},\do,x_s)=(x_{e+u+1},\do,x_s).$$
Then 
$$U(x_*)=U(x'_*)\T U(s_*)\T U(x''_*),\qua 
U(y_*)=U(x'_*)\T U(s'_*)\T U(x''_*).$$
 Now $\tU(x_*)$ (resp. $\tU(y_*)$) is the set of orbits for the action of a
subgroup of $U_A\T U_A$ on $\tU(x'_*)\T\tU(s_*)\T\tU(x''_*)$ (resp. $\tU(x'_*)\T\tU(s'_*)\T\tU(x''_*)$) given by 
$$(u,u'):(\x,\x',\x'')\m((1,u)\x,(u,u')\x',(u',1)\x'')$$
 (the subgroup is $U_A\T U_A$ if $x'_*,x''_*$ are nonempty,
is $\{1\}\T U_A$ if $x'_*$ is empty and $x''_*$ is nonempty, is $U_A\T\{1\}$ if $x'_*$ is nonempty and $x''_*$ is
empty, is $\{1\}\T\{1\}$ if $x'_*$ and $x''_*$ are empty.) Let $\ct$ be
the set of orbits of the same subgroup of $U_A\T U_A$ on $\tU(x'_*)\T(U_A\dw U_A)\T\tU(x''_*)$ for the action 
given by the same formulas as above. We have a diagram
$$U(x_*)@>h>>U(x'_*)\T(U^w_A\dw)\T U(x''_*)@<h'<<U(y_*)$$
where 
$$h(\x,\x',\x'')=(\x,\ph_{s_*}(\x'),\x''),\qua 
h'(\x,\x',\x'')=(\x,\ph_{s'_*}(\x'),\x'')$$
(with $\ph_{s_*}$,
$\ph_{s'_*}$ as in 2.7(g). Note that $h,h'$ are bijections. We set $H=h'{}\i h$. We have a diagram
$$\tU(x_*)@>\ti h>>\ct@<\ti h'<<\tU(y_*)$$
where 
$$\ti h(\x,\x',\x'')=(\x,\ps_{s_*}(\x'),\x''),\qua
\ti h'(\x,\x',\x'')=(\x,\ps_{s'_*}(\x'),\x'')$$ 
(with $\ps_{s_*}$, $\ps_{s'_*}$ as in 2.7(i)). Note that $\ti h,\ti h'$ are bijections (they are well defined
by 2.7(i)). We set $\tH=(\ti h')\i\ti h$. It is clear that $H,\tH$ satisfy the requirements of (a). This 
proves (a).

\subhead 2.10\endsubhead
Let $\d$ be an automorphism of $\car$ that is, a triple consisting of automorphisms $\d:Y@>>>Y$, $\d:X@>>>X$ and a
bijection $\d:I@>>>I$ such that $\la\d(y),\d(x)\ra=\la y,x\ra$ for $y\in Y,x\in X$, $\d$ is compatible with the 
imbeddings $I@>>>Y$, $I@>>>X$ and $\d(i)\cdot\d(j)=i\cdot j$ for $i,j\in I$. There is a unique group automorphism
of $W$, $w\m\d(w)$ such that $\d(s_i)=s_{\d(i)}$ for all $i\in I$. There is a unique algebra automorphism 
(preserving $1$) of $\bof_A$, $x\m\d(x)$, such that $\d(\th_i^{(c)})=\th_{\d(i)}^{(c)}$ for all 
$i\in I,c\in\NN$. We have $\d(\BB)=\BB$. There is a unique algebra automorphism of $\dUU_A$, $u\m\d(u)$, such 
that $\d(b^-1_\z b'{}^+)=\d(b)^-1_{\d(\z)}\d(b')^+$ for all $b\in\BB,b'\in\BB,\z\in X$. We have $\d(\dBB)=\dBB$.
There is a unique algebra automorphism (preserving $1$) of $\hUU_A$, $u\m\d(u)$ such that  
$\d(\sum_{a\in\dBB}n_aa)=\sum_{a\in\dBB}n_a\d(a)$ for all functions $\dBB@>>>A$, $a\m n_a$. This automorphism 
restricts to a group automorphism $G_A@>>>G_A$ denoted again by $\d$ and to an automorphism of $U_A$. For any
$i\in I,h\in H$ we have $\d(x_i(h))=x_{\d(i)}(h)$, $\d(\ds_i)=\ds_{\d(i)}$.
 
\subhead 2.11\endsubhead
Let $A\in\cc,A'\in\cc$ and let $\c:A@>>>A'$ be a homomorphism of rings preserving $1$. There is a unique 
ring homomorphism (preserving $1$) $\hUU_A@>>>\hUU_{A'}$, $u\m\c(u)$ such that
$\c(\sum_{a\in\dBB}n_aa)=\sum_{a\in\dBB}\c(n_a)a$ for all functions $\dBB@>>>A$, $a\m n_a$. This restricts to
a group homomorphism $G_A@>>>G_{A'}$ denoted again by $\c$ and to a group homomorphism $U_A@>>>U_{A'}$. For any 
$i\in I,h\in A$ we have $\c(x_i(h))=x_i(\c(h))$, $\c(\ds_i)=\ds_i$.

\head 3. The main results\endhead
\subhead 3.1\endsubhead
In this section $A\in\cc$ is fixed unless otherwise specified. We also fix an automorphism $\d$ of $\car$ as in 
2.10 and a ring automorphism $\c$ of $A$ preserving $1$. There are induced group automorphisms of $G_A$ denoted 
again by $\d,\c$ (see 2.10, 2.11). These automorphisms commute; we set $\p=\d\c=\c\d:G_A@>>>G_A$; note that $\p$ 
maps $U_A$ onto itself.

Two elements $w,w'$ of $W$ are said to be $\d$-conjugate if $w'=y\i w\d(y)$ for some $y\in W$. The relation of
$\d$-conjugacy is an equivalence relation on $W$; the equivalence classes are said to be $\d$-conjugacy classes. 
A $\d$-conjugacy class $C$ in $W$ (or an element of it) is said to be $\d$-elliptic if $C\cap W_J=\em$ for any
$J\sneq I,\d(J)=J$. Let $C$ be a $\d$-elliptic $\d$-conjugacy class in $W$. Let $C_{min}$ be the set of elements 
of minimal length of $C$. 

For any $w\in W$ we define a map 
$$\X_A^w:U_A\T(U^w_A\dw)@>>>U_A\dw U_A$$ 
by $(u,z)\m uz\p(u)\i$.

\subhead 3.2\endsubhead
In this subsection we assume that $x,y\in W$ are such that $l(x\d(y))=l(x)+l(y)=l(yx)$. We show: 

($*$) {\it $\X^{x\d(y)}_A$ is injective if and only if $\X^{yx}_A$ is injective.}
\nl
Let $x_*=(x,\d(y)),x'_*=(y,x)$.
In the following proof we write $U,U^x,U^y,U^{\d(y)}$ instead of $U_A,U^x_A,U^y_A,U^{\d(y)}_A$. 
Now $\X^{x\d(y)}_A$ can be identified with 
$$U\T(U^x\dx)\T(U^{\d(y)}\d(\dy))@>>>\tU(x_*),\qua (u,z,z')\m[uz,z'\p(u)\i]$$
and $\X^{yx}_A$ can be identified with 
$$U\T(U^y\dy)\T(U^x\dx)@>>>\tU(x'_*),\qua (u,z',z)\m[uz',z\p(u)\i].$$
(We use 2.7(g), 2.7(i).) The condition that $\X^{x\d(y)}_A$ is injective can be stated as follows: 

(a) {\it if $u\in U,z_1,z_2\in U^x\dx,z_3,z_4\in U^{\d(y)}\d(\dy)$ satisfy $uz_1v\i=z_2,vz_3\p(u)\i=z_4$ for
some $v\in U$ then $u=1$.}
\nl
The condition that $\X^{yx}_A$ is injective can be stated as follows: 

(b) {\it if $u'\in U,z'_1,z'_2\in U^x\dx,z'_3,z'_4\in U^y\dy$ satisfy 
$u'z'_3v'{}\i=z'_4, v'z'_1\p(u')\i=z'_2$ for some $v'\in U$ then $u'=1$.}
\nl
Assume that (a) holds and that the hypothesis of (b) holds. We have $v'z'_1\p(u')\i=z'_2$, 
$\p(u')\p(z'_3)\p(v')\i=\p(z'_4)$ and $\p(z'_3),\p(z'_4)\in U^{\d(y)}\d(\dy)$. Applying (a) with $v=\p(u')$, 
$u=v'$ we obtain $v'=1$. Then $\p(u')=z'_2{}\i z'_1\in\dx\i U^x\dx\sub U^-$. But $U\cap U^-=\{1\}$ (see 2.3(b)) 
hence $u'=1$. Thus the conclusion of (b) holds. 

Next we assume that (b) holds and that the hypothesis of (a) holds. We have $\p\i(v)\p\i(z_3)u\i=\p\i(z_4)$ and 
$\p\i(z_3),\p\i(z_4)\in U^y\dy$. Applying (b) with $v'=u$, $u'=\p\i(v)$ we obtain $v=1$. Then
$\p(u)=z_4\i z_3\in\dy\i U^y\dy\sub U^-$. But $U\cap U^-=\{1\}$ hence $u=1$. Thus the conclusion of 
(a) holds. We see that (a) holds if and only if (b) holds. This proves ($*$).

\subhead 3.3\endsubhead
In this subsection we assume that $x,y\in W$ are such that $l(yx)=l(x)+l(y)$ and we write $U,U^x,U^y,U^{yx}$ 
instead of $U_A,U^x_A,U^y_A,U^{yx}_A$. We have a commutative diagram
$$\CD
U\T U\T(U^y\dy)\T(U^x\dx)@>\Th_A>>(U\dy U)\T(U\dx U)\\
@VmVV     @Vm'VV\\
U\T(U^{yx}\dy\dx)@>\X^{yx}_A>>U\dy\dx U\endCD$$
where 
$$\Th_A(u,u',g,g')=(ugu'{}\i,u'g'\p(u)\i),\qua m(u,u',g,g')=(u,gg'),\qua m'(g,g')=gg'.$$
We show 

(a) {\it the sets 
$$S'=\{\X':U\dy\dx U@>>>U\T(U^{yx}\dy\dx);\X_A^{yx}\X'=1\},$$
$$S''=\{Z:(U\dy U)\T(U\dx U)@>>>U\T U\T(U^y\dy)\T(U^x\dx);\Th_AZ=1\}$$
are in natural bijection.}
\nl
Let $\X'\in S'$. We define $Z$ as follows. Let $(z,z')\in(U\dy U)\T(U\dx U)$. We have $\X'(zz')=(u,gg')$ where 
$u\in U$, $g\in U^y\dy$, $g'\in U^x\dx$ are uniquely defined. Since $ug\in U\dy U$ we can write uniquely 
$ug=z_0u'$ where $z_0\in U^y\dy$, $u'\in U$ (see 2.7(d)). We set $Z(z,z')=(u,u',g,g')$. This defines the map $Z$.
We have 
$$zz'=\X_A^{yx}(\X'(zz'))=ugg'\p(u)\i.$$
The equality $zz'=(ug)(g'\p(u)\i)$ and 2.7(h) imply that
$ugv\i=z,vg'\p(u)\i=z'$ for a well defined $v\in U$. We have $ug=zv=z_0u'$. From $zv=z_0u'$ and 2.7(d) we see that
$u'=v$. Thus $ugu'{}\i=z,u'g'\p(u)\i=z'$. Hence 
$$\Th_A(Z(z,z'))=\Th_A(u,u',g,g')=(ugu'{}\i,u'g'\p(u)\i)=(z,z').$$ 
Thus $\Th_AZ=1$ and $Z\in S''$.

Conversely, let $Z\in S''$. We define $\X'$ as follows. Let $h\in U\dy\dx U$. Using 2.7(h) and 2.7(d) we can 
write uniquely $h=zz'$ where $z\in U^y\dy$, $z'\in U\dx U$. We set $\X'(h)=m(Z(z,z'))$. We have 
$$\X^{yx}_A(\X'(h))=\X^{yx}_A(m(Z(z,z'))=m'(\Th_A(Z(z,z'))=m'(z,z')=zz'=h.$$ 
Thus $\X^{yx}_A\X'=1$ and $\X'\in S'$.

It is easy to check that the maps $S'@>>>S''$, $S''@>>>S'$ defined above are inverse to each other. This proves
(a).

\subhead 3.4\endsubhead
In this subsection we assume that $x,y\in W$ are such that 
$$l(x\d(y))=l(x)+l(y)=l(yx)$$
and we write $U,U^x,U^y,U^{\d(y)},U^{yx},U^{x\d(y)}$ instead of \lb 
$U_A,U^x_A,U^y_A,U^{\d(y)}_A,U^{yx}_A,U^{x\d(y)}_A$. We show 

(a) {\it the sets 
$$S'_1=\{\X'_1:U\dy\dx U@>>>U\T(U^{yx}\dy\dx);\X_A^{yx}\X'_1=1\},$$
$$S'_2=\{\X'_2:U\dx\d(\dy)U@>>>U\T(U^{x\d(y)}\dx\d(\dy));\X_A^{yx}\X'_2=1\}$$
are in natural bijection.}
\nl
Define $\Th:U\T U\T(U^y\dy)\T(U^x\dx)@>>>(U\dy U)\T(U\dx U)$ and

$\Th':U\T U\T(U^x\dx)\T(U^{\d(y)}\d(\dy))@>>>(U\dx U)\T(U\d(\dy)U)$
\nl
by $(u,u',g,g')\m(ugu'{}\i,u'g'\p(u)\i)$. In view of 3.3 applied to $x,y$ and also to $\d(y),x$ we see that 
to prove (a), it is enough to show:

(b) {\it the sets 
$$S''_1=\{Z:(U\dy U)\T(U\dx U)@>>>U\T U\T(U^y\dy)\T(U^x\dx);\Th Z=1\},$$
$$S''_2=\{Z':(U\dx U)\T(U\d(\dy)U)@>>>U\T U\T(U^x\dx)\T(U^{\d(y)}\d(\dy));\Th'Z'=1\}$$
are in natural bijection.}
\nl
here 

Now (b) follows from the commutative diagram
$$\CD
U\T U\T(U^y\dy)\T(U^x\dx)@>\Th>>(U\dy U)\T(U\dx U)\\
@V\io VV    @V\io'VV   \\
U\T U\T(U^x\dx)\T(U^{\d(y)}\d(\dy))@>\Th'>>(U\dx U)\T(U\d(\dy)U)\endCD$$
where $\io(u,u',g,g')=(u',\p(u),g',\p(g))$ and $\io'(g,g')=(g',\p(g))$ are bijections. This proves (a).

\subhead 3.5\endsubhead
We show:

(a) {\it if $\X^w_A$ is injective for some $w\in C_{min}$ then $\X^w_A$ is injective for any $w\in C_{min}$;}
\nl
For any $w,w'$ in $C_{min}$ there exists a sequence $w=w_1,w_2,\do,w_r=w'$ in $C_{min}$ such that for any 
$h\in[1,r-1]$ we have either $w_h=x\d(y)$, $w_{h+1}=yx$ for some $x,y$ as in 3.2 or $w_{h+1}=x\d(y)$, $w_h=yx$ 
for some $x,y$ as in 3.2 (See \cite{\GP}, \cite{\GKP}, \cite{\HE}.) Now (a) follows by applying 3.2($*$) several 
times.

Now assume that for some $w'\in C_{min}$, we are given $\X':U_A\dw'U_A@>>>U_A\T(U^{w'}_A\dw')$ such that 
$\X_A^{w'}\X'=1$. Let $w\in C_{min}$. We show how to
construct a map $\X'':U_A\dw U_A@>>>U_A\T(U^w_A\dw)$ such that $\X_A^w\X''=1$. We choose a sequence
$w=w_1,w_2,\do,w_r=w'$ in $C_{min}$ as in the proof of (a). We define a sequence of maps
$\X'_i:U_A\dw_i U_A@>>>U_A\T(U^{w_i}_A\dw_i)$ ($i\in[1,r]$) such that $\X_A^{w_i}\X'_i=1$ by induction on $i$
as follows. We set $\X'_1=\X'$. Assuming that $\X'_i$ is defined for some $i\in[1,r-1]$ we define $\X'_{i+1}$ so 
that $\X'_i,\X'_{i+1}$ correspond to each other under a bijection as in 3.4. Then the map $\X'':=\X'_r$ satisfies
our requirement. In particular, we see that:

(b) {\it if $\X^w_A$ is surjective for some $w\in C_{min}$ then $\X^w_A$ is surjective for any $w\in C_{min}$.}

\proclaim{Theorem 3.6} Recall that $A\in\cc$. Let $C$ be a $\d$-elliptic $\d$-conjugacy class in $W$ and
let $w\in C_{min}$. Then:

(i) $\X_A^w$ is injective;

(ii) if $\c=1$, then $\X_A^w$ is bijective;

(iii) if $A$ is a field, $\c$ has finite order $m$ and the fixed point field $A^\c$ is perfect, then $\X_A^w$ is
bijective;

(iv) if $A$ is an algebraic closure of a finite field $F_q$ and $\c(x)=x^q$ for all $x\in A$ then $\X_A^w$ is 
bijective;

(v) if $A$ is finite and $\c$ is arbitrary then $\X_A^w$ is bijective.
\endproclaim
The proof will occupy 3.7-3.10. If $W=\{1\}$ the result is trivial. Hence we can assume that
$W\ne\{1\}$. 

\subhead 3.7\endsubhead
Let $e$ be the order of any element of $C$. According to \cite{\GM}, \cite{\GKP}, \cite{\HE}, we can find an 
element $w'\in C_{min}$ such that $\hw'\widehat{\d(w')}\do\widehat{\d^{e-1}(w')}=\hy_1\hy_2\do\hy_t$
(in the braid monoid $\hW$) where $y_*=(y_1,y_2,\do,y_t)$ is a sequence in $W$ such that $y_1=w_I$.
Since $W\ne\{1\}$ we have $t\ge2$. From now until the end of 3.9 we assume that $w=w'$.
Let $x_*=(w,\d(w),\do,\d^{e-1}(w))$ and let $H:U(x_*)@>\si>>U(y_*)$, 
$\tH:U(x_*)@>\si>>U(y_*)$ be as in 2.9(a). Let $z\in U_A\dw U_A$. Then $[z,\p(z),\do,\p^{e-1}(z)]\in\tU(x_*)$ 
hence $\tH[z,\p(z),\do,\p^{e-1}(z)]\in\tU(y_*)$. We set $u=\z(\tH[z,\p(z),\do,\p^{e-1}(z)])\in U_A$ where
$\z:U(y_*)@>>>U$ is as in 2.8. We can write uniquely $\p^{-e}(u)z\p^{-e+1}(u)\i=z'u'$ where $z'\in U_A^w\dw$, 
$u'\in U_A$ (see 2.7(d)). We set $\X'_A(z)=(\p^{-e}(u)\i,z')\in U_A\T(U^w_A\dw)$. Thus we have a map 
$$\X'_A:U_A\dw U_A@>>>U_A\T(U^w_A\dw).$$
We show:

(a) {\it $\X'_A\X_A^w$ is the identity map of $U_A\T(U^w_A\dw)$ into itself.}
\nl
Let $(u_1,z_1)\in U_A\T(U^w_A\dw)$. Let $z=u_1z_1\p(u_1)\i$. We must show that $\X'_A(z)=(u_1,z_1)$. Let 
$\x=\tH[z,\p(z),\do,\p^{e-1}(z)]\in\tU(y_*)$. By 2.8(b) we have $\x=(1,u\i)[h_1,\do,h_t]$ with 
$(h_1,\do,h_t)\in U(y_*)$, $u\in U_A$ uniquely determined; moreover from the definitions we have $u=\z(\x)$.
We have
$$\align&[z,\p(z),\do,\p^{e-1}(z)]\\&
=[u_1z_1\p(u_1)\i,\p(u_1)\p(z_1)\p^2(u_1)\i,\do,\p^{e-1}(u_1)\p^{e-1}(z_1)\p^e(u_1)\i]\\&
=[u_1z_1,\p(z_1),\do, \p^{e-2}(z_1),\p^{e-1}(z_1)\p^e(u_1)\i]\\&
=(u_1,\p^e(u_1))[z_1,\p(z_1),\do,\p^{e-2}(z_1),\p^{e-1}(z_1)]\in\tU(x_*),\endalign$$
hence 
$$\align&\x=\tH((u_1,\p^e(u_1))[z_1,\p(z_1),\do, \p^{e-2}(z_1),\p^{e-1}(z_1)])\\&=
(u_1,\p^e(u_1))\tH([z_1,\p(z_1),\do, \p^{e-2}(z_1),\p^{e-1}(z_1)])\\&=
(u_1,\p^e(u_1))[a_1,\do,a_t]\endalign$$ 
where 
$(a_1,\do,a_t)=H(z_1,\p(z_1),\do, \p^{e-2}(z_1),\p^{e-1}(z_1))\in U(y_*)$.
Thus we have 
$$(1,u\i)[h_1,\do,h_t]=(u_1,\p^e(u_1))[a_1,\do,a_t]\in\tU(y_*)$$
 and
$$[h_1,\do,h_t]=[u_1a_1,a_2,\do,a_{t-1},a_t\p^e(u_1)\i u\i]\in\tU(y_*).$$
Note that $h_i\in U^{y_i}_A\dy_i$ for $i\in[1,t]$, $a_i\in U^{y_i}_A\dy_i$ for $i\in[2,t]$ and
$u_1a_1\in U_A^{y_1}\dy_1$ (we use that $y_1=w_I$, $U_A^{w_I}=U_A$). Using the uniqueness part of
2.8(b) we deduce that $\p^e(u_1)\i u\i=1$, hence $\p^{-e}(u)=u_1\i$. Then we have 
$$\p^{-e}(u)z\p^{-e+1}(u)\i=u_1\i z\p(u_1)=z_1\in U_A^w\dw.$$
Using the definition we have $\X'_A(z)=(\p^{-e}(u)\i,z_1)=(u_1,z_1)$. This proves (a).

\subhead 3.8\endsubhead
If $w$ is as in 3.7 then from 3.7(a) we see that $\X_A^w$ is injective. This proves 3.6(i) for this $w$.

\subhead 3.9\endsubhead
Let $w$ be as in 3.7. We set $N=l(w)+l(w_I)$.
We identify $A^N=U_A\T(U^w_A\dw)$ as in 2.6(a), 2.7(e) and $A^N=U_A\dw U_A$ as in 2.7(f). 
Then $\Xi_A^w$ and $\Xi'_A$ become maps  $f_A:A^N@>>>A^N$, $f'_A:A^N@>>>A^N$ such that
$f'_Af_A=1$. Assuming that $\c=1$ and $\d$ is fixed, we see from the definitions that $(f_A)_{A\in\cc}$, 
$(f'_A)_{A\in\cc}$ are polynomial families, see 1.1. (Note that the definition of $f'_A$ involves the 
isomorphisms $H,\ti H$ in 2.8(a). By the proof of 2.8(a) these isomorphisms can be regarded as
polynomial families when $A$ varies.) We can now apply 1.3 and we see 
that $f_A$ is bijective for any $A$. This proves 3.6(ii) for our $w$.

Now assume that $A,\c,m,A^\c$ are as in 3.6(iii). Let $A_0$ be an algebraic closure of $A_1:=A^\c$. Let 
$A_2=A\ot_{A_1}A_0\in\cc$. 
Now $f_A:A^N@>>>A^N$ is not given by polynomials with coefficients in $A$; however, $A$ is an $A_1$-vector space 
of dimension $m$ and $f_A$ can be viewed as a map $A_1^{Nm}@>>>A_1^{Nm}$ given by polynomials with coefficients
in $A_1$. The same polynomials describe the map $f_{A_2}:A_2^N@>>>A_2^N$ viewed as a map $A_0^{Nm}@>>>A_0^{Nm}$.
This last map is injective by 3.8 (applied to $A_2$) and then it is automatically bijective by \cite{\BR} (see 
also \cite{\AX}, \cite{\GRR, 10.4.11}) applied to the affine space $A_0^{Nm}$ over $A_0$. Thus $f_{A_2}$ is 
bijective. Let $\x\in A^N$.
Then $\x':=f_{A_2}\i(\x)\in A_2^N$ is well defined. The Galois group of $A_0$ over $A_1$ acts on $A_2$ (via the 
action on the second factor) hence on $A_2^N$. This action is compatible with $f_{A_2}$ and it fixes $\x$ hence 
it fixes $\x'$. Since $A_1$ is perfect it follows that $\x'\in(A\ot_{A_1}A_1)^N=A^N$. We have $f_A(\x')=\x$. Thus
$f_A$ is surjective, hence bijective. This proves 3.6(iii) for our $w$.

Now assume that $A,\c$ are as in 3.6(iv). In this case $f_A$ can be viewed as a map $A^N@>>>A^N$ given by 
polynomials with coefficients in $A$. This map is injective by 3.8 and then it is automatically bijective by
\cite{\BR} (see also \cite{\AX}, \cite{\GRR, 10.4.11}) applied to the affine space $A^N$. Thus $f_A$ is 
bijective. This proves 3.6(iv) for our $w$.

Finally assume that $A$ is finite. Then $A^N$ is finite. Since $f_A$ is injective, it is automatically bijective.
This proves 3.6(v) for our $w$.

\subhead 3.10\endsubhead
Now let $w$ be any element of $C_{min}$. 

Since $\Xi^{w'}_A$ is injective for $w'$ as in 3.7 (see 3.8) we see using 3.5(a) that $\Xi^w_A$ is injective.
This proves 3.6(i) for our $w$.

Now assume that we are in the setup of 3.6(ii),(iii),(iv) or (v). Since $\X^{w'}_A$ is bijective for $w'$ as in 
3.7 (see 3.9) we see using 3.5(a),(b) that $\X^w_A$ is bijective. This completes the proof of Theorem 3.6.

\subhead 3.11\endsubhead
Recall that $A\in\cc$. 
In this subsection we assume that $\c=1$. Let $w\in C_{min}$ ($C$ as in 3.6). 
We identify $A^N=U_A\T(U^w_A\dw)$ (with $N=l(w)+l(w_I)$)
as in 2.6(a), 2.7(e) and $A^N=U_A\dw U_A$ as in 2.7(f). Then $\X_A^w$ becomes a map $A^N@>>>A^N$. From the
definitions we see that $(\X_{A'}^w)_{A'\in\cc}$ is a polynomial family. (Here $\d$ is fixed and $\c=1$ for any 
$A'$.) Hence $(\X_A^w)^*:A[X_1,\do,X_N]@>>>A[X_1,\do,X_N]$ is defined for any $A\in\cc$. We have the following
result.

\proclaim{Theorem 3.12} In the setup of 3.11, $(\X_A^w)^*$ is an isomorphism. In particular, if $A$ is an
algebraically closed field, then $\X_A^w:A^N@>>>A^N$ is an isomorphism of algebraic varieties.
\endproclaim
Let $w'\in C_{min}$ be as in 3.7. As we have seen in 3.9, there exists a polynomial family $\X'_A:A^N@>>>A^N$ 
($A\in\cc$) such that $\X'_A\X^{w'}_A=1$ for any $A$. Since $\X^{w'}_A$ is bijective by 3.6(ii) we must also have
$\X^{w'}_A\X'_A=1$. Now the method in the paragraph preceding 3.5(b) yields for any $A$ an explicit map 
$\X''_A:A^N@>>>A^N$ such that $\X^w_A\X''_A=1$. Moreover from the definitions we see that $(\X''_A)_{A\in\cc}$ is
a polynomial family. It follows that $(\X''_A)^*$ is defined and $(\X''_A)^*(\X^w_A)^*=1$. Since $\X^w_A$ is a 
bijection (see 3.6(ii)), we see that $\X^w_A\X''_A=1$ implies $\X''_A\X^w_A=1$ hence $(\X^w_A)^*(\X''_A)^*=1$. 
This, together with $(\X''_A)^*(\X^w_A)^*=1$ shows that $(\X^w_A)^*$ is an isomorphism. This completes the 
proof of the theorem. Note that the second assertion of the theorem can alternatively be proved using 1.2(ii) and
the fact that $\X_A^w$ is injective (see 3.6(i)).

\subhead 3.13\endsubhead
In this subsection we assume that $A,\c$ are as in 3.6(iv). Let $w\in C_{min}$ ($C$ as in 3.6). We identify 
$A^N=U_A\T(U^w_A\dw)$ (with $N=l(w)+l(w_I)$)
as in 2.6(a), 2.7(e) and $A^N=U_A\dw U_A$ as in 2.7(f). Then $\X_A^w$ becomes a map 
$A^N@>>>A^N$. It is in fact a morphism of algebraic varieties. Exactly as in 3.12 we define a map 
$\X''_A:A^N@>>>A^N$ such that $\X^w_A\X''_A=1$. But this time $\X''_A$ is not a morphism of algebraic varieties 
but only a quasi-morphism (see \cite{\XW, 1.1}). (This is because the definition of $\X''_A$ involves 
$\c\i:A@>>>A$ which is a quasi-morphism but not a morphism.) Since $\X^w_A$ is a bijection (see 3.6(iv)) we 
deduce that we have also $\X''_A\X^w_A=1$. Thus we have the following result:

(a) {\it $\X^w_A:A^N@>>>A^N$ is a bijective morphism whose inverse is a quasi-morphism.}

\subhead 3.14\endsubhead
Let $C$ be as in 3.6 and let $w\in C_{min}$. Let $A\in\cc$ and let $\d,\c,\p$ be as in 2.1. We define

(a) {\it $\a^w_A:{}^{\d\i(w)\i}U_A\T U^w_A@>>>U_A$ by $(u',u'')\m u'u''\dw\p(u')\i\dw\i$.}
\nl
We show:

(b) {\it $\a^w_A$ is injective.}
\nl
Assume that 

$(u',u'')\in{}^{\d\i(w)\i}U_A\T U^w_A$, $(u'_1,u''_1)\in{}^{\d\i(w)\i}U_A\T U^w_A$
\nl
and $u'u''\dw\p(u')\i\dw\i=u'_1u''_1\dw\p(u'_1)\i\dw\i$. Then $\X^w_A(u',u''\dw)=\X^w_A(u'_1,u''_1\dw)$. Using 
3.6(i) we deduce $(u',u''\dw)=(u'_1,u''_1\dw)$ hence $u'=u'_1$, $u''=u''_1$ and (a) follows.

We show:

(c) {\it If $\c,A$ are as in 3.6(ii),(iii),(iv) or (v), then $\a^w_A$ is bijective.}
\nl
Let $u\in U_A$. By 3.6, we have $u\dw=u'u''\dw\p(u')\i$ for some $u'\in U_A,u''\in U^w_A$. By 2.7(c) we can write
$u'{}\i u=u_1u_2$ where $u_1\in U^w_A,u_2\in{}^wU_A$. Then $u_1\dw(\dw\i u_2\dw)=u''\dw\p(u')\i$. Using 2.7(d) we
deduce that $\p(u')\i=\dw\i u_2\dw\in{}^{w\i}U_A$. Thus $u'\in{}^{\d\i(w)\i}U_A$ so that $\a^w_A$ is surjective. 
Together with (b) this implies that $\a^w_A$ is bijective.

We note:

(d) {\it If $\c=1$ and $A$ is an algebraically closed field then $\a^w_A$ is an isomorphism of algebraic 
varieties.}
\nl
We note that $(\a^w_{A'})_{A'\in\cc}$ can be viewed as a polynomial family of injective maps $A'{}^n@>>>A'{}^n$ 
($A'\in\cc$, $n$ as in 2.1). Hence the result follows from 1.2(ii).

\head 4. Applications\endhead
\subhead 4.1\endsubhead
In this section we assume that $A$ is an algebraically closed field. We write $G,U,{}^wU,U^w,T$ instead of 
$G_A,U_A,{}^wU_A,U^w_A,T_A$. By \cite{\LU, 4.11}, $G$ is naturally a connected reductive 
algebraic group over $A$ with root datum $\car$ and $U$ is the unipotent radical of a Borel subgroup $B^*$ of 
$G$ with maximal torus $T$ normalized by each $\ds_i$. We assume that $G$ is semisimple or equivalenty that
$\{i';i\in I\}$ span a subgroup of finite index in $X$. Let $\d$ be an automorphism of $\car$ (necessarily of
finite order, say $c$). The corresponding group automorphism $\d:G@>>>G$ (see 2.10) preserves the algebraic 
group structure and has finite order $c$. Let $\hG$ be the semidirect product of $G$
with the cyclic group of order $c$ with generator $d$ such that $dxd\i=\d(x)$ for all $x\in G$. 
Then $\hG$ is an algebraic group with identity component $G$. Let $\cb$ be the variety of Borel subgroups
of $G$. For each $w\in W$ let $\co_w$ be the set of all $(B,B')\in\cb\T\cb$ such that $B=xB^*x\i$,
$B'=x\dw B^*\dw\i x\i$ for some $x\in G$. We have $\cb\T\cb=\sqc_{w\in W}\co_w$. As in \cite{\XW, 0.1, 0.2}, for 
any $w\in W$ let

$\fB_w=\{(g,B)\in Gd\T\cb;(B,gBg\i)\in\co_w\}$,

$\tfB_w=\{(g,g'{}^wU)\in Gd\T G/{}^wU;g'{}\i gg'\in\dw Ud\}$.
\nl
Define $\p_w:\tfB_w@>>>\fB_w$ by $(g,g'{}^wU)\m(g,g'B^*g'{}\i)$.

In the remainder of this section we assume that $C$ is a $\d$-elliptic $\d$-conjugay class in $W$ and that
$w\in C_{min}$. Then $\p_w$ is a principal bundle with group $T_w=\{t_1\in T;\dw\i t\dw=dtd\i\}$, a finite 
abelian group (see {\it loc.cit.}); the group $T_w$ acts on $\tfB_w$ by $t:(g,g'{}^wU)\m(g,g't\i{}^wU)$. Now 
$G$ acts on $\fB_w$ by $x:(g,B)\m(xgx\i,xBx\i)$ and on $\tfB_w$ by $x:(g,g'{}^wU)\m(xgx\i,xg'{}^wU)$. We show:

(a) {\it Let $\co$ be a $G$-orbit in $\tfB_w$. There is a unique $v\in U^{\d(w)}$ such that \lb
$(\dw vd,{}^wU)\in\co$.}
\nl
Clearly $\co$ contains an element of the form $(\dw ud,{}^wU)$  where $u\in U$. We first show the existence of
$v$. It is enough to show that for some $z\in{}^wU,v\in U^{\d(w)}$ we have $z\dw udz\i=\dw vd$ that is
$u= \dw\i z\i\dw v\d(z)$. Setting $z'=\dw\i z\i\dw$, $w'=\d(w)$ we see that it is enough to show that
$u= z'v\dw'\d(z')\i\dw'{}\i$ for some $z'\in {}^{\d\i(w')\i}U, v\in U^{w'}$. But this follows from 3.14(c)
with $\c=1$ and $w$ replaced by $w'$.

Now assume that $(\dw vd,{}^wU)\in\co$, $(\dw v'd,{}^wU)\in\co$ where $v,v'\in U^{\d(w)}$. We have 
$\dw v'd=u\dw vdu\i$ for some $u\in{}^wU$. Setting $u'=\dw\i u\dw\in{}^{w\i}U$ we have $v'd\dw=u' vd\dw u'{}\i$, 
that is $v'=u'v\d(\dw)\d(u'{}\i)\d(\dw)\i$. Using 3.14(b) (appplied to $\d(w)$ instead of $w$) we see that $u'=1$
and $v=v'$. This completes the proof of (a).

We can reformulate (a) as follows.

(b) {\it The closed subvariety $\{(\dw vd,{}^wU);v\in U^{\d(w)}\}$ of $\tfB_w$ meets each $G$-orbit in $\tfB_w$ in
exactly one point. Hence the space of $G$-orbits in $\tfB_w$ can be identified with the affine space $U^{\d(w)}$.}
\nl
We show:

(c) {\it The closed subvariety $\{(\dw vd,B^*); v\in U^{\d(w)}\}$ of $\fB_w$ is isomorphic to $U^{\d(w)}$; its
intersection with any $G$-orbit in $\fB_w$ is a single $T_w$-orbit (for the restriction of the $G$-action), hence
is a finite nonempty set.}
\nl
The first assertion of (c) is obvious. Now let $\bco$ be a $G$-orbit in $\fB_w$. Let 
$$Z=\{(\dw vd,B^*);v\in U^{\d(w)}\}\cap\bco.$$ 
There exists a $G$-orbit $\co$ in $\tfB_w$ such that 
$\p_w(\co)=\bco$. By (b) we can find $v\in U^{\d(w)}$ such that $(\dw vd,{}^wU)\in\co$. Then
$(\dw vd,B^*)=\p_w(\dw vd,{}^wU)\in\bco$ so that $Z\ne\em$. If $(\dw vd,B^*)\in Z$ and $t\in T_w$ then 
$(t\dw vdt\i,B^*)\in\bco$ and $(t\dw vdt\i,B^*)=(\dw v''d,{}^wU)$ where 
$$v''=\dw\i t\dw v'dt\i d\i=(dtd\i)v'(dt\i d\i)\in{}^wU.$$
Thus $(t\dw vdt\i,B^*)\in Z$ so that $T_w$ acts on $Z$.

Assume now that $v,v'\in U^{\d(w)}$ are such that $(\dw vd,B^*)\in\bco$, $(\dw v'd,B^*)\in\bco$. Then for some 
$x\in G$ we have 
$$\p_w(x\dw vdx\i,x{}^wU)=\p_w(\dw v'd,{}^wU).$$
Since $\p_w$ is a principal fibration with group $T_w$ it follows that 
$$(x\dw vdx\i,x{}^wU)=(\dw v'd,t\i{}^wU)$$ 
for some $t\in T_w$. Thus
$(\dw vd,{}^wU),(t\dw v'dt\i,{}^wU)$ are in the same $G$-orbit on $\tfB_w$. Note that 
$(t\dw v'dt\i,{}^wU)=(\dw v''d,{}^wU)$ where 
$$v''=\dw\i t\dw v'dt\i d\i=(dtd\i)v'(dt\i d\i)\in{}^wU.$$
Using (b) we deduce that $v=v''$. Thus 
$$(\dw v'd,B^*)=(\dw(\dw\i t\i\dw)v(dtd\i)d,B^*)=(t\i\dw vdt,B^*)=t\i(\dw vd,B^*)$$
so that $(\dw v'd,B^*),(\dw vd,B^*)$ are in the same $T_w$-orbit. This completes the proof of (c).

We can reformulate (c) as follows.

(d) {\it The closed subvariety $\{(\dw vd,B^*);v\in U^{\d(w)}\}$ of $\fB_w$ meets each $G$-orbit in $\fB_w$ in 
exactly one $T_w$-orbit. Hence the space of $G$-orbits in $\fB_w$ can be identified with the orbit space of the 
affine space $U^{\d(w)}$ under an action of the finite group $T_w$.}
\nl
Statements like the last sentence in (b) and (d) were proved in \cite{\XW, 0.4(a)} assuming that $G$ is almost 
simple of type $A,B,C$ or $D$. The extension to exceptional types is new.

\subhead 4.2\endsubhead
In this subsection we assume that $\d=1$ so that $d=1$. Let $\g$ be the unipotent class of $G$ attached to $C$ in
\cite{\WEU}. Recall from {\it loc.cit.} that $\g$ has codimension $l(w)$ in $G$. The following result exhibits a 
closed subvariety of $G$ isomorphic to the affine space $A^{l(w)}$ which intersects $\g$ in a finite set.

(a) {\it The closed subvariety $\Si:=\dw U^w$ of $G$ is isomorphic to $U^w$ and $\Si\cap\g$ is a 
single $T_w$-orbit (for the conjugation action), hence is a finite nonempty set.}
\nl
According to \cite{\HO}, the subset $\fB^\g_w=\{(g,B)\in\fB_w;g\in\g\}$ is a single $G$-orbit on $\fB_w$. Let 
$Z=\{(\dw v,B^*); v\in U^w\}\cap\fB^\g_w$, $Z'=\{\dw v; v\in U^w\}\cap\g$. The first projection defines a
surjective map $Z@>>>Z'$. Since $Z$ is a single $T_w$-orbit (see 3.1(c)), it follows that $Z'$ is a single 
$T_w$-orbit. This proves (a).

\subhead 4.3\endsubhead
In this subsection we assume that $A,F_q,\c$ are as in 3.6(iv). Then $\p=\d\c:G@>>>G$ is the Frobenius map for an 
$F_q$-rational structure on $G$. As in \cite{\DL}, we set

$\tX_w=\{g'{}^wU\in G/{}^wU;g'{}\i\p(g')\in\dw U\}$.
\nl
Now the finite group $G^\p:=\{g\in G;\p(g)=g\}$ acts on $\tX_w$ by by $x:g'{}^wU\m xg'{}^wU$. Let
${}^wU\bsl\bsl U$ be the set of orbits of the ${}^wU$-action on $U$ given by $u_1:u\m\dw\i u_1\dw u\p(u_1)\i$.
According to \cite{\DL, 1.12}, we have a bijection $G^\p\bsl\tX_w@>\si>>{}^wU\bsl\bsl U$, 
$g'{}^wU\m\dw\i g'{}\i\p(g')$ with inverse induced by $u\m g'{}^wU$ where $g'\in G$, $g'{}\i\p(g')=\dw u$. Under 
the substitution $\dw\i u_1\dw=u'$, the ${}^wU$-action above on $U$ becomes the ${}^{w\i}U$-action on $U$ given by
$u_2:u\m u' u\d(\dw)p(u')\i\d(\dw)\i$. Using 3.14(c) for $\d(w)$ instead of $w$ we see that the space of orbits 
of this action can be identified with $U^{\d(w)}$. Thus we have the following result.

(a) {\it The space of orbits of $G^\p$ on $\tX_w$ is quasi-isomorphic to the affine space $U^{\d(w)}$.}
\nl
A statement like (a) was proved in \cite{\XW} assuming that $G$ is almost simple of type $A,B,C$ or $D$
and $\d=1$. The extension to general $G$ is new.


\widestnumber\key{GKP}  
\Refs
\ref\key{\AX}\by J.Ax\paper Injective endomorphisms of algebraic varieties and schemes\jour Pacific J.Math.\vol31
\yr1969\pages1-7\endref
\ref\key{\BR}\by A.Bialynicki-Birula and M.Rosenlicht\paper Proc.Amer.Math.Soc.\vol13\yr1962\pages200-203\endref
\ref\key{\DL}\by\by P.Deligne and G.Lusztig\paper Representations of reductive groups over finite fields\lb\jour 
Ann.Math.\vol103\yr1976\pages103-161\endref
\ref\key{\GP}\by M.Geck and G.Pfeiffer\book Characters of finite Coxeter groups and representations of \lb 
Iwahori-Hecke algebras\bookinfo LMS Monographs\vol21\yr2000\publ Oxford Univ.Press\endref
\ref\key{\GKP}\by M.Geck, S.Kim and G.Pfeiffer\paper Minimal length elements in twisted conjugacy classes of 
finite Coxeter groups\jour J.Algebra\vol229\yr2000\pages570-600\endref
\ref\key{\GM}\by M.Geck and J.Michel\paper Good elements of finite Coxeter groups and representations of 
Iwahori-Hecke algebras\jour Proc.London Math.Soc.\vol74\yr1997\pages275-305\endref
\ref\key{\GRR}\by A.Grothendieck\paper \'El\'ements de g\'eom\'etrie alg\'ebrique IV, \'Etude locale des sch\'emas et des morphismes de sch\'emas (Troisi\`eme partie)\jour Inst.Hautes \'Et.Sci. Publ.Math.\vol28\yr1966\endref
\ref\key{\GR}\by A.Grothendieck\paper \'El\'ements de g\'eom\'etrie alg\'ebrique IV, \'Etude locale des sch\'emas
et des morphismes de sch\'emas (Quatri\`eme partie)\jour Inst.Hautes \'Et.Sci. Publ.Math.\vol32\yr1967\endref
\ref\key{\HE}\by X.He\paper Minimal length elements in some double cosets of Coxeter groups\jour Adv.Math.\vol215
\yr2007\pages469-503\endref
\ref\key{\QG}\by G.Lusztig\book Introduction to quantum groups\bookinfo Progr.in Math.110\publ Birkh\"auser 
Boston\yr1993\endref 
\ref\key{\LU}\by G.Lusztig\paper Study of a $\ZZ$-form of the coordinate ring of a reductive group\jour 
Jour. Amer. Math. Soc.\vol22\yr2009\pages739-769\endref
\ref\key{\WEU}\by G.Lusztig\paper From conjugacy classes in the Weyl group to unipotent classes\jour
arxiv:1003.0412\endref 
\ref\key{\HO}\by G.Lusztig\paper Elliptic elements in a Weyl group: a homogeneity property\jour arxiv:1007.5040
\endref
\ref\key{\XW}\by G.Lusztig\paper On certain varieties attached to a Weyl group element\jour arxiv:1012.2074\endref
\ref\key{\SE}\by A.Sevostyanov\paper Algebraic group analogues of the Slodowy slices and deformations of Poisson 
$W$-algebras\paperinfo doi:10.1093/imrn/rnq139 \jour Int. Math. Res. Notices\yr2010\endref
\ref\key{\ST}\by R.Steinberg\paper Regular elements of semisimple algebraic groups\jour Inst. Hautes \'Et.Sci. 
Publ. Math.\vol25\yr1965\pages49-80\endref
\endRefs
\enddocument